\documentclass[11pt]{article}
\usepackage{amssymb,amsmath,amsthm}
\newtheorem{theorem}{Theorem}
\newtheorem{cor}[theorem]{Corollary}
\newtheorem{prop}[theorem]{Proposition}
\newtheorem{lemma}[theorem]{Lemma}
\newcommand{\C}{\mathcal{C}}

\newcommand{\R}{\,\overline{\! R}}
\begin{document}

\title{{\huge \bf Asymmetric binary covering codes}}
\author{{\Large Joshua N.\ Cooper}\\
{\it Department of Mathematics, University of California } \\ {\it
at San Diego, La Jolla, California}\\ E-mail:
jcooper@math.ucsd.edu\\ \\
{\Large Robert B.\ Ellis}\\
{\it Department of Mathematics, Texas A\&M University}\\
E-mail: rellis@math.tamu.edu\\
\\
and\\ \\
{\Large Andrew B.\ Kahng}\\
{\it Department of Computer Science and Engineering,
University of California  } \\ {\it at San Diego, La Jolla, California}\\
E-mail: abk@cs.ucsd.edu}

\date{May 3, 2002}
\maketitle

\begin{abstract}
An asymmetric binary covering code of length $n$ and radius $R$ is
a subset $\mathcal{C}$ of the $n$-cube $Q_n$ such that every
vector $x\in Q_n$ can be obtained from some vector
$c\in\mathcal{C}$ by changing at most $R$ 1's of $c$ to 0's, where
$R$ is as small as possible. $K^+(n,R)$ is defined as the smallest
size of such a code. We show $K^+(n,R)\in \Theta(2^n/n^R)$ for
constant $R$, using an asymmetric sphere-covering bound and
probabilistic methods. We show $K^+(n,n-\R)=\R+1$ for constant
coradius $\R$ iff $n\geq \R(\R+1)/2$. These two results are
extended to near-constant $R$ and $\R$, respectively. Various
bounds on $K^+$ are given in terms of the total number of 0's or
1's in a minimal code.  The dimension of a minimal asymmetric
linear binary code ($[n,R]^+$-code) is determined to be
$\min{\!\{0,n-R\}}$. We conclude by discussing open problems and
techniques to compute explicit values for $K^+$, giving a table of
best known bounds.
\end{abstract}

\section{Introduction} Suppose we wish to have a small set of binary
$n$-vectors with the property that {\it every} binary $n$-vector
is no more than $R$ bit flips from one of them.  This is the
classical question of finding ``covering codes.''  Recent surveys
of results on covering codes appear in \cite{BLP98} and
\cite{CHLL97}, and earlier important results appear in
\cite{CKMS85,GS85}.  The topic of covering codes continues to be
an active area of research, and the interested reader is referred
to \cite{L02} for a comprehensive bibliography of the subject.

Let $Q_n$ be the set of binary $n$-vectors
$\{x=(x_1,x_2,\ldots,x_n):x_i\in \{0,1\}\}$ with algebraic
structure inherited from the vector space $\mathbb{F}_2^n$ and the
partial ordering inherited from the boolean lattice (i.e.,
$x\preceq y$ if $x_i\leq y_i$ for all $1\leq i\leq n$). We denote
the ``top'' and ``bottom'' elements, i.e., $(1,\ldots,1)$ and
$(0,\ldots,0)$, by $\hat{1}$ and $\hat{0}$, respectively. Define
the {\em weight}, or {\em level}, of $x \in Q_n$ as
$w(x)=\sum_{i=1}^n x_i$, where each coordinate is treated as an
ordinary integer (equivalently, $w(x)$ is the number of ones in
$x$). Define the Hamming distance between $x$ and $y$ as
$d(x,y)=w(x+y)$.  The \textit{undirected ball} in $Q_n$ with
center $x$ and radius $R$, denoted by $B_n(x,R)$, is the set
$\{y\in Q_n:d(x,y)\leq R\}$. The {\em covering radius} of a set
$\mathcal{C}\subseteq Q_n$ is the smallest integer $R\geq 0$ such
that $Q_n=\cup_{c\in \mathcal{C}}B_n(c,R)$. The ordinary
definition of a binary covering code, which for our purposes we
refer to as a {\em symmetric binary covering code\/} of length $n$
and radius $R$, or more simply an $(n,R)$-code, is a set of
``codewords'' $\mathcal{C}\subseteq Q_n$ with covering radius $R$.
We use $K(n,R)$ to denote the minimum size of any $(n,R)$-code.

We now consider the additional restriction of requiring the bit
flips used to go from a vector to its covering codeword to be in
only one direction. This restriction arises in a problem of layout
data compression in VLSI design which motivated the present work
\cite{EKZ02}.  Data encoding the placement of certain metal
features on a microchip can be transmitted with at most $R$ errors
per $n$ bits using a covering code, except that metal may only be
removed (to an extent controlled by $R$) and not added, so as to
avoid causing a short-circuit. This simple variation on ordinary
(``symmetric'') covering codes opens up a world of questions, with
many of the answers quite different from the symmetric case.  The
dual problem of ``unidirectional'' error-correcting/detecting
codes has been studied in \cite{R74,RSSS94,SSA96}.

The extra restriction is now formalized in the definition of an
asymmetric covering code. The \textit{upward directed ball} in
$Q_n$ with center $x$ and radius $R$ is $B_n^+(x,R)=
B_n(x,R)\cap\{y\in Q_n:x\preceq y\}$, and the corresponding
\textit{downward directed ball} is $B_n^-(x,R)=B_n(x,R)\cap \{y\in
Q_n:y\preceq x\}$.  We write $b_n^+(x,R)$ and $b_n^-(x,R)$ for the
sizes of the directed balls $B_n^+(x,R)$ and $B_n^-(x,R)$,
respectively.  We sometimes instead say $b_n^+(l,R)$ or
$b_n^-(l,R)$, where $l$ is the weight of $x$, since $b_n^+$ and
$b_n^-$ depend only on the weight of the ball's center.  Indeed,
we have
\begin{equation}\label{eqn:ballSizeBound}
b^+_n(l,R) = b^-_n(n-l,R) = \sum_{j=0}^R \binom{n-l}{j} \leq
\sum_{j=0}^R \binom{n-l}{j}\binom{R}{R-j} =  \binom{n-l+R}{R}.
\end{equation}

A set $\C\subseteq Q_n$ {\em downward $R$-covers} $Q_n$ provided
that $Q_n=\cup_{c\in\C}B^-_n(c,R)$, and the {\em asymmetric
covering radius} of $\C$ is the smallest $R$ for which $\C$
downward $R$-covers $Q_n$. We define an {\em asymmetric binary
covering code} of length $n$ and radius $R$, or more simply an
$(n,R)^+$-code, to be any set $\C\subseteq Q_n$ with covering
radius $R$. We sometimes refer to the \textit{coradius} $\R:=n-R$
of an $(n,R)^+$-code when $R$ is large. Our main object of study
is the function $K^+(n,R)$, defined to be the minimum size of any
$(n,R)^+$-code.

Finally, denote the concatenation of two vectors
$x=(x_1,\ldots,x_n)$ and $y=(y_1,\ldots,y_m)$ by $(x|y) =
(x_1,\ldots,x_n,y_1,\ldots,y_m)$, in precisely the same way as it
is defined for symmetric codes.  The {\em direct sum} of two sets
$X$ and $Y$ is $X\oplus Y:=\{(x|y):x\in X, y\in Y\}$. Note that if
$\C$ is an $(n_1,R_1)^+$-code and $\C^\prime$ is an
$(n_2,R_2)^+$-code, then $\C \oplus \C^\prime$ is an
$(n_1+n_2,R_1+R_2)^+$-code, and so
\begin{equation}\label{eqn:directSumBound}
K^+(n,R) \leq \ K^+(n_1,R_1)\cdot K^+(n-n_1,R-R_1).
\end{equation}
We will use this observation several times in the course of our
discussion.

In this paper, we explain several substantive differences and
similarities between symmetric and asymmetric binary covering
codes, and offer directions for further investigation. Section
\ref{sec:asymptotics} gives the exact asymptotic order of
magnitude of the size of minimal codes with constant radius and
gives exact asymptotics in the case of constant coradius. The
bounds we provide are then used to derive somewhat weaker bounds
in a completely general setting. The topic of Section
\ref{sec:difference} is the increase that the size of a minimal
code experiences when its length or radius is incremented or
decremented, respectively.  We tackle linear asymmetric codes in
Section \ref{sec:linear} -- a surprisingly simple matter, given
the complexity of the issue in the symmetric case -- and we finish
with several open problems and a table of our best known bounds in
Section \ref{sec:conclusion}.

\section{Asymptotic bounds} \label{sec:asymptotics}

We can achieve a lower bound for the asymptotic order of magnitude
of $K^+(n,R)$ for constant $R$ by considering a variant of the
traditional sphere-covering bound. Sphere-covering lower bounds
are achieved by examining the size of the (directed or undirected)
balls of a given radius centered at each vector. The
straightforward sphere-covering bound in the symmetric case
appears as \cite[Theorem 6.1.2]{CHLL97}, which we state here for
completeness, and then extend to the asymmetric case.
\begin{theorem}[Sphere-covering bound]
\begin{equation}
K(n,R) \geq
\left\lceil\frac{2^n}{\sum_{j=0}^R\binom{n}{j}}\right\rceil.
\label{eqn:sphereCoveringBound}
\end{equation}
\end{theorem}

\begin{theorem}[Asymmetric sphere-covering bound]
\label{thm:asymmetricSphereCoveringLowerBound} Let $0\leq R\leq
n$.  Then
\begin{equation} K^+(n,R) \geq \left\lceil \sum_{l=0}^n
\frac{\binom{n}{l}}{\sum_{j=0}^R\binom{\min(n,l+R)}{j}}\right\rceil.
\label{eqn:fractionalLowerBound}
\end{equation}
\end{theorem}
\proof For any $(n,R)^+$-code $\C$, we may write
$$
|\C| = \sum_{c \in \C} 1 = \sum_{c \in \C} \sum_{v \in B^-_n(c,R)}
    b^-_n(w(c),R)^{-1}.
$$
Switching the order of summation yields
$$
|\C| = \sum_{v \in Q_n} \sum_{c \in B^+_n(v,R) \cap \C} b^-_n(w(c),R)^{-1}
$$
For a vector $v$ of weight $l$, the
largest directed ball of radius $R$ that could contain it is
centered at a vector of weight $l+R$ and has size
$\sum_{j=0}^R\binom{l+R}{j}.$ However, if $l+R > n$, then the
largest ball that could contain $v$ is the one which is centered
at $\hat{1}$.  Therefore, since every vector in $Q_n$ must be covered
by at least one $c \in \C$,
\begin{align*}
|\C| & \geq \sum_{v \in Q_n} \sum_{c \in B^+_n(v,R) \cap \C} \left(
\sum_{j=0}^R \binom{\min(n,w(v)+R)}{j} \right )^{-1} \\
& \geq \sum_{v \in Q_n} \left (\sum_{j=0}^R \binom{\min(n,w(v)+R)}{j}
\right )^{-1}.
\end{align*}
Noting that $\binom{n}{l}$ vertices have weight $l$ gives the desired
result.
\qed

The desired lower bound for $K^+(n,R)$ is determined by bounding
the denominator of each term in (\ref{eqn:fractionalLowerBound}).
Using the bound on ball size from (\ref{eqn:ballSizeBound}), for
all $0 \leq l \leq n$ and a fixed $R$,
$$
b^-_n(l,R) \leq \binom{n+R}{R} \in O\!\left( n^R \right).
$$
Therefore, we find the following as an immediate consequence of
Theorem \ref{thm:asymmetricSphereCoveringLowerBound}.

\begin{cor}\label{cor:asymmetricSphereCoveringAsymptoticLowerBound} Fix $R
\geq 0$.  Then
$$
K^+(n,R) \in \Omega \! \left( \frac{2^n}{n^R} \right).
$$
\end{cor}

We can count more carefully than we did in Theorem
\ref{thm:asymmetricSphereCoveringLowerBound} by specifying a
system of inequalities that the code must satisfy. For an
arbitrary $(n,R)^+$-code, define the sequence
$(a_0,a_1,\ldots,a_n)$ by letting $a_l$ be the number of codewords
of weight $l$. Whenever necessary, define
$a_{n+1}=\ldots=a_{n+R}=0$.  We have the following lemma.

\begin{lemma}
Let $\C$ be an $(n,R)^+$-code. Then the number of codewords $a_l$
of weight $l$ must satisfy
\begin{equation}
    a_l \geq \binom{n}{l}-\sum_{j=1}^Ra_{l+j}\binom{l+j}{j}.
    \label{eqn:modifiedSphereBoundLemma}
\end{equation}
\end{lemma}
\proof There are $\binom{n}{l}$ vertices on level $l$ to be
covered.  At most $a_{l+j}\binom{l+j}{j}$ of these points can be
$R$-covered by the $a_{l+j}$ codewords of weight $l+j$. The rest
must be included as codewords themselves. \qed

In fact, (\ref{eqn:modifiedSphereBoundLemma}) must hold for all
$0\leq l\leq n$ for all $(n,R)^+$-codes.  Therefore we may
construct an integer program based on these restrictions to
provide another lower bound for $K^+(n,R)$. This is what was done
(with some minor refinements) to find most of the lower bounds
presented in the table at the end of this paper.

\begin{prop}\label{prop:modifiedOneSideSphereCoverBound}
Let $IP^+(n,R)$ be the result of the following integer program.
\begin{eqnarray}
\mathrm{Minimize}\qquad \sum_{l=0}^n a_l   && \mathrm{subject \
to}\nonumber \\
\sum_{j=0}^{R}a_{l+j}\binom{l+j}{j} & \geq & \binom{n}{l},
\qquad\mathrm{for}\mbox{\, $0\leq l\leq n$},
\label{eqn:modifiedSphereBoundConstraint}\\
\mathrm{and}\qquad a_l & \geq & 0, \ \ \mathrm{integer,}
\qquad\mathrm{for}\mbox{\, $0\leq l\leq n$}. \nonumber
\end{eqnarray}
Then $K^+(n,R)\geq IP^+(n,R)$.  Furthermore, this bound is at
least as good as the asymmetric sphere-covering bound
(\ref{eqn:fractionalLowerBound}).
\end{prop}
\proof We have already established that the program yields a lower
bound for $K^+(n,R)$.  It remains to show that it is at least as
good as (\ref{eqn:fractionalLowerBound}).  Let $(a_0, \ldots,
a_n)$ be a solution vector which achieves $IP^+(n,R)$.  Applying
the conventions that $b_n^-(l,R)=b_n^-(n,R)$ when $l\geq n$, all
indices vary over those integers not excluded explicitly, and
$a_l=0$ for $l<0$ and $l>n$, we have
\begin{align*}
\sum_l a_l & = \sum_l \sum_{i \leq R} \frac{a_l \binom{l}{i}}{b^-_n(l,R)}
\\
& = \sum_{l^\prime} \sum_{i \leq R} \frac{a_{l^\prime + i} \binom{l^\prime
+ i}{i}}{b^-_n(l^\prime + i,R)}
\end{align*}
by making the substitution $l^\prime = l - i$.  Then, using the
fact that downward ball sizes are monotone increasing in the
weight of their centers, and applying
(\ref{eqn:modifiedSphereBoundConstraint}),
\begin{align*}
\sum_l a_l & \geq \sum_{l^\prime} \sum_{i \leq R} \frac{a_{l^\prime + i}
\binom{l^\prime + i}{i}}{b^-_n(l^\prime + R,R)} \\
& \geq \sum_{l^\prime} \frac{\binom{n}{l^\prime}}{b^-_n(l^\prime + R,R)}
\end{align*}
which is the asymmetric sphere covering bound. \qed

Note that the IP is relatively small, since its coefficient matrix is just
$(n+1) \times (n+1)$.

\subsection{Asymptotic order of magnitude for small radius} \label{subsec:constRadius}

Clearly, $K^+(n,0)=2^n$ since all 0-balls contain only their
centers. For positive $R$, however, the issue is much more
complicated.  In particular, we wish to understand the growth of
$K^+(n,R)$ in $n$ for constant $R$.  The lower bound given by
Corollary \ref{cor:asymmetricSphereCoveringAsymptoticLowerBound}
says that the density $|\C|/2^n$ of a minimal $(n,R)^+$-code $\C$
is $\Omega(n^{-R})$; the probabilistic arguments given in this
section show that this is, in fact, achievable.

Define a \textit{patched asymmetric covering code of radius $R$},
or a patched $(n,R)^+$-code, to be a pair $(S,T)$ with $S,T
\subset Q_n$ such that $S$ has covering radius $R$ with respect to
covering only $Q_n\setminus T$.  Thus every vector in the cube is
either in $B_n^-(s,R)$ for some $s\in S$ or in the ``patch'' $T$.
We say that the $\delta$-weight of the patched asymmetric covering
code $(S,T)$ is $|S|+\delta|T|$, and define $p(n,R,\delta)$ to be
the minimum $\delta$-weight over all patched $(n,R)^+$-codes.

For a given patched $(n,R)^+$-code $(S,T)$ and a $(k,R)^+$-code
$\C$, we define the \textit{semi-direct sum} of $(S,T)$ and $\C$,
denoted by $(S,T) \boxplus \C$, to be $\left(S\oplus Q_k\right)
\cup \left(T\oplus \C\right)$. It is easy to verify that $(S,T)
\boxplus \C$ is an $(n+k,R)^+$-code.  In the next two
propositions, we will generate a small cover for $Q_{2n}$ first by
showing in Proposition \ref{prop:deletionCover} the existence of a
patched $(n,R)^+$-code $(S,T)$ with low $\delta$-weight, and then
by building in Proposition \ref{prop:inductiveCover} a
$(2n,R)^+$-code $\C$ from the semi-direct sum of the patched
$(n,R)^+$-code and a small $(n,R)^+$-code found inductively.

\begin{prop} \label{prop:deletionCover} Let $R\geq 0$ be fixed.  For some
absolute constant
$\alpha_R > 0$ and any $\delta > 0$,
$$ p(n,R,\delta) \leq \frac{\alpha_R
2^n}{n^R} \big(\max\{\log(\delta n^R / \alpha_R),0\} + 1\big).
$$
\end{prop}
\proof A standard argument (pointed out by J. Bell \cite{B01})
permits us to choose $\alpha_R$ to be the least real number so
that for all positive integers $n$,
$$
\sum_{j=0}^{n} \binom{n}{j} b_n^+(j,R)^{-1}  \leq \alpha_R
\frac{2^n}{n^R}.
$$
If $\delta < \alpha_R / n^R$, then choosing $T$ be be all of $Q_n$
yields the desired result.  Thus, we may assume $\delta \geq
\alpha_R / n^R$.  Randomly choose a patched asymmetric cover
$(S,T)$ as follows. Let $p_j = \min\{\log(\delta n^R / \alpha_R)
b_n^+(j,R)^{-1},1\}$ for $j=0 \ldots n-1$, and let $p_n = 1$. For
each vector $v$ in the cube, add it to $S$ with probability
$p_{w(v)}$.  Then, add all the uncovered points to $T$.  The
expected $\delta$-weight of $(S,T)$ is, by linearity of
expectation,
$$
\mathbf{E}(|S|) + \delta \mathbf{E}(|T|) \leq \sum_{j=0}^n
\binom{n}{j} p_j + \delta \sum_{v \in Q_n} \mathbf{P}(v
\mathrm{\mbox{ is uncovered}}).
$$
The probability that a vector $v$ is uncovered is the product,
over each of the $b_n^+(v,R)$ vertices that could cover $v$, of
the probability that each vertex is not chosen. Thus,
$$
\mathbf{P}(v \mathrm{\mbox{ is uncovered}}) \ = \prod_{u\in
B^+_n(v,R)}(1-p_{w(u)}) \ \leq \ (1-p_{w(v)})^{b^+_n(v,R)}
$$
and we have (using the formula $(1-1/x)^x\leq e^{-1}$ which is
valid for all $x\geq 1$)
\begin{align*}
\mathbf{E}(|S|) + \delta \mathbf{E}(|T|) &\leq
    \log\left(\frac{\delta n^R}{\alpha_R}\right) \sum_{j=0}^{n}
    \frac{\binom{n}{j}}{b^+_n(j,R)}  + \delta
    \sum_{j=0}^{n} \binom{n}{j} (1-p_j)^{b^+_n(j,R)}\qquad \\
& \leq \log\left(\frac{\delta n^R}{\alpha_R}\right) \alpha_R
\frac{2^n}{n^R} + \delta \sum_{j=0}^{n}
\binom{n}{j} \frac{\alpha_R}{\delta n^R} \\
\qquad\qquad\qquad\qquad\qquad & \leq \alpha_R \frac{2^n}{n^R}
\left(\log\left(\frac{\delta n^R}{\alpha_R}\right) + 1\right),
\end{align*}
and so there exists a patched cover of the desired
$\delta$-weight. \qed

This leads immediately to the following.

\begin{prop} \label{prop:inductiveCover}
For each $R \geq 0$, there exists a $\beta_R > 0$ such that for
every nonnegative integer $m$,
$$
K^+\left(2^m,R\right) \leq \frac{\beta_R 2^{2^m}}{2^{mR}}.
$$
\end{prop}
\proof If $R=0$, the result is trivial, so we may assume that
$R>0$. We proceed by induction on $m$. We require the constant
$\beta_R$ defined as
$$ \beta_R = \max\left\{ \frac{1}{2},
    \min\left\{x:x\geq \alpha_R \mbox{ and } x\geq
    2^R\alpha_R\left(\log(x/\alpha_R)+1\right)
    \right\} \right\}.
$$
The statement certainly holds for $m =0$, since $K^+(1,R) = 1$.
Assume it is true for $m$.  We construct a cover for $Q_{2^{m+1}}$
by taking the semi-direct sum of a patched $(2^m,R)^+$-code
$(S,T)$ achieving $p\left(2^m,R,K^+(2^m,R)/2^{2^m}\right)$ and a
minimal $(2^m,R)^+$-code. The result is a $(2^{m+1},R)^+$-code of
size
$$
|(S,T) \boxplus \C| = |S|\cdot 2^{2^m}+|T|\cdot K^+(2^m,R) =
    2^{2^m} p\left(2^m,R,\frac{K^+(2^m,R)}{2^{2^m}}\right).
$$
By the previous proposition, then,
$$ K^+(2^{m+1},R) \leq 2^{2^m} \alpha_R \frac{2^{2^m}}{2^{mR}}
    \left(\max\left\{\log\left(\frac{K^+(2^m,R)2^{mR}}
    {\alpha_R 2^{2^m}}\right),0\right\} + 1\right).
$$
If we apply the inductive hypothesis to bound $K^+(2^m,R)$, we
find that
$$
K^+(2^{m+1},R) \leq
    \frac{2^{2^{m+1}}}{(2^{m+1})^R} 2^R \alpha_R
    \left(\max\left\{\log\left(\frac{\beta_R}
        {\alpha_R}\right),0\right\}
    + 1\right) \leq \beta_R \frac{2^{2^{m+1}}}{(2^{m+1})^R}
$$
by the choice of $\beta_R$.\qed

A straightforward application of the direct sum formula
(\ref{eqn:directSumBound}) allows us to generalize this result to
all nonnegative integers $n$ from those which are powers of 2.

\begin{cor} Let $R \geq 0$ be fixed.  For some absolute constant
$\gamma_R > 0$ and every integer $n$,
$$
K^+(n,R) \leq \frac{\gamma_R 2^n}{n^R}.
$$
\end{cor}
\proof Set $\gamma_R = 2^R \beta_R$, and let $m = \lfloor
\log_2(n) \rfloor$.  Then by the direct sum formula,
\begin{eqnarray*}
K^+(n,R) & \leq & K^+(2^m,R) \cdot K^+(n-2^m,0) = K^+(2^m,R)
    \cdot 2^{n-2^m}     \nonumber \\
& \leq &  \frac{\beta_R 2^n}{2^{mR}} \leq
    \frac{2^R \beta_R 2^n}{n^R} =
    \frac{\gamma_R 2^n} {n^R}, \nonumber
\end{eqnarray*}
since $n/2^m\leq 2$. \qed

This, combined with Corollary
\ref{cor:asymmetricSphereCoveringAsymptoticLowerBound}, gives us
the following characterization of the asymptotic behavior of
$K^+(n,R)$ for $R$ constant.

\begin{theorem} For a fixed $R \geq 0$, $K^+(n,R) \in \theta(2^n/n^R)$.
\end{theorem}

The probabilistic technique used above, reminiscent of the
so-called ``deletion method'' (see, for instance \cite{AS00}),
applies in a more general setting, although the results are
significantly weaker. The following proof is essentially a very
simple version of the proof of Proposition
\ref{prop:deletionCover}, with $\delta=1$, but we include it
because the bound achieved is of independent interest. Define
$\nu(n,R)$ by
$$
\nu(n,R) = \sum_{j=0}^n \frac{\binom{n}{j}}{b^+_n(j,R)}.
$$
Then we have the following, analogous to \cite[Theorem
12.1.2]{CHLL97}.

\begin{prop} For any $n,R \geq 0$, $K^+(n,R) \leq (n \log 2 + 1)
\nu(n,R)$.
\end{prop}
\proof We construct an $(n,R)^+$-code probabilistically.  Add
vectors $x \in Q_n$ to $\C$ independently with probability
$$
p_{w(x)} =
\min\left\{1,\frac{\log(2^n/\nu(n,R))}{b^+_n(w(x),R)}\right\},
$$
and then add to $\C$ all points not covered by directed balls of
radius $R$ centered at the chosen $x$'s.  Just as in the proof of
Proposition \ref{prop:deletionCover}, the resulting code has
expected size
\begin{align*}
\mathbf{E}(|\C|) & \leq \sum_{j=0}^n \binom{n}{j} p_j + \sum_{j=0}^n
\binom{n}{j} (1-p_j)^{b^+_n(j,R)} \\
& \leq \left(\log(2^n/\nu(n,R)) + 1\right) \,
    \nu(n,R) \leq (n \log 2 + 1) \, \nu(n,R).
\end{align*}
Therefore, there exists an $(n,R)^+$-code of at most this size.
\qed

The preceding proposition can be used to achieve upper bounds in
specific cases comparable to those known for symmetric covering
codes. For example, a routine calculation gives that $K^+(n,R)$ is
within a multiplicative factor of $O(n)$ of the asymmetric
sphere-covering lower bound whenever $R \in O(\sqrt{n})$.

\subsection{Asymptotics for large radius}

Since $R \geq n$ implies $K^+(n,R) = 1$, another region of
interest in the space of possible $n$'s and $R$'s is the case of
constant (and positive) coradius $\R = n-R$.  The very precise
asymptotics we achieve for this case also permit a much rougher
analysis of $K^+(n,R)$ for general $n$ and $R$, with the best
results occurring when $R$ is  close to $n$.

A few trivial values are immediate. $K^+(n,n)=1$, since the
downward $n$-ball at $\hat{1}$ covers everything, and
$K^+(n,n-1)=2$ by considering the code
$\C=\{(1,1,\ldots,1),(0,1,\ldots,1)\}$.  In fact, for fixed $\R$,
the sequence $\{K^+(n,n-\R)\}_n$ converges to $\R+1$ in a manner
we now characterize.

\begin{lemma}\label{lem:diagCover}
$K^+(n,n-\R)\leq \R+1$ for \ $n\geq \R(\R+1)/2$ \ and $\R\geq 0$.
\end{lemma}
\proof Construct the $(n,n-\R)^+$-code
$$\C =
\{(1,1,1,1,\ldots,1),(0,1,1,1,\ldots,1),(1,0,0,1,\ldots,1),\ldots\}$$
of size $\R+1$, where the $(i+1)^{\mbox{{\small{st}}}}$ codeword
has $i$ consecutive 0's starting in position $(i-1)i/2+1$. Having
$n\geq \R(\R+1)/2$ is required in order for there to be enough
positions to place all the 0's. To see that $\C$ is an
$(n,n-\R)^+$-code, let $x\in Q_n$. If $w(x)\geq \R$, $x$ is
covered by $(1,\ldots,1)$. Otherwise, $x$ could only avoid being
covered by the $w(x)+1$ codewords on levels
$(n-\R+w(x)),\ldots,(n-\R)$ by having, for each codeword $c \in
\C$ on these levels, a 1 in a position where $c$ has a 0.  This is
impossible, since $x$ has $w(x)$ 1's and the positions of 0's in
the $w(x)+1$ codewords are disjoint. \qed

The same set $\C$ is a symmetric $(n,n-\R)$-code, and therefore
yields an upper bound on $K(n,n-\R)$ for $n\geq \R(\R+1)/2$. This
bound is not tight in the symmetric case, but is in fact tight in
the asymmetric case, due to the additional structure imposed by
requiring
 vectors to be covered by codewords above them.  (In fact, $K(n,R)=2$
 whenever $R+1 \leq n \leq 2R+1 $.)  We
summarize the behavior of $K^+(n,n-\R)$ in the theorem below.
Before we proceed, however, we have the following definition and
lemma.

For a set of vectors $\C \subset Q_n$ and an index $i \in
\{1,\ldots,n\}$, define $\C^i \subset Q_{n-1}$, the
\textit{contraction of $\C$ at i}, to be the set of points of $\C$
with a ``1'' at position $i$, projected into $Q_{n-1}$ by deletion
of that bit.

\begin{lemma} \label{lem:RestrictCodes} If $\C$ downward
$R$-covers $Q_n$, then for each $i \in \{1,\ldots,n\}$, $\C^i$
downward $R$-covers $Q_{n-1}$.
\end{lemma}
\proof Let $Q^\prime_n$ be the set of points in $Q_n$ with a ``1''
at coordinate $i$.  Note that if $x \in Q^\prime_n$, $y \in Q_n$,
and $y \succ x$, then $y \in Q^\prime_n$.  Thus the union of all
the downward-directed Hamming balls of radius $R$ centered at the
points of $\C \cap Q^\prime_n$ must contain $Q^\prime_n$.  If we
project $Q^\prime_n$ onto $Q_{n-1}$ in the natural way, then the
image of $\C \cap Q^\prime_n$ downward $R$-covers $Q_{n-1}$. \qed

\bigskip

\begin{theorem} \label{thm:superDiag} $K^+(n,n-\R)\geq \R+1$ for $n\geq 1$
and $\R\geq 0$, with
equality when $n\geq n_{\R}:=\R(\R+1)/2$.   Furthermore, $n_{\R}$
is the least integer $n$ for which equality holds.
\end{theorem}
\proof  First, we show by induction on $\R$ that if $\C$ downward
$(n-\R)$-covers $Q_n$, then there are at least $x$ codewords of
$\C$ with weight $< n-\R+x$, for all nonnegative $x\leq \R$. This
is certainly true for $\R=0$; assume it is true for $\R-1$. Now
consider a general $\R>0$. There is at least one codeword $c$ at
or below level $n-\R$, since we must cover the vertex $\hat{0}$.
Because $\R>0$, we may choose some coordinate $i$ where $c$ has a
zero. Lemma \ref{lem:RestrictCodes} gives that the contraction
$\C^i$ downward $((n-1)-(\R-1))$-covers $Q_{n-1}$. By induction,
there are at least $x$ points in $\C^i$ with weight $<
(n-1)-(\R-1)+x=n-\R+x$, for all $x\leq \R-1$.  However, level $l$
in the $(n-1)$-cube corresponds to level $l+1$ in the original
cube, which has an additional codeword at or below level $n-\R$.
Therefore, $\C$ has at least $x$ codewords with weight $<n-\R+x$
for each $x$ with $0 \leq x \leq \R$.  Taking $x=\R$ gives us the
desired lower bound, and combining this with Lemma
\ref{lem:diagCover} yields $K^+(n,n-\R)=\R+1$ for $n \geq n_{\R}$.
It remains to show that $K^+(n,n-\R)>\R+1$ for $n<n_{\R}$.

To that end, suppose $n<n_{\R}$.  Since the number of codewords
below level $n-\R+x$ (or, having at least $(\R-x)$ 0's) is at
least $x$, a minimal code $\C$ achieving $K^+(n,n-\R)$ has at
least $\sum_{x=0}^{\R} (\R-x)=n_{\R}$ total 0's in its codewords.
Since $n<n_{\R}$, two of the codewords must have 0's in a common
position. If we contract $\C$ at that coordinate, the resulting
code is at least 2 smaller than the original one, and so
$K^+(n,n-\R) \geq K^+(n-1,n-\R) + 2 \geq (\R-1)+1 + 2 = \R+2$.
\qed

By combining this theorem with the direct sum construction, we can
bound $K^+(n,R)$ from above for a wide range of parameters.

\begin{theorem} For any nonnegative $n$ and $\R$, $K^+(n,n-\R) \leq
 \left (2n/\R \right )^{\lceil \R^2/(2n-\R) \rceil}$.
\end{theorem}
\proof Note that by the direct sum construction, for any $0 \leq n
\leq n^\prime$ and $0 \leq R \leq R^\prime$, we have
$$
K^+(n,R) \leq K^+(n,R) \cdot K^+(n^\prime - n, R^\prime - R) \leq
K^+(n^\prime,R^\prime),
$$
so that $K^+(n,R)$ is nondecreasing in both parameters. Therefore,
applying the direct sum construction again for any integer $M >
0$,
\begin{align*}
K^+(n,n-\R) &\leq K^+\left (M \lceil n/M \rceil, M (\lceil n/M \rceil -
\lfloor \R/M \rfloor) \right ) \\
&\leq K^+\left (\lceil n/M \rceil, \lceil n/M \rceil - \lfloor
\R/M \rfloor \right )^M.
\end{align*}
If we choose $M \geq \R^2/(2n-\R)$, then it is straightforward to
see that
$$
\left \lceil \frac{n}{M} \right \rceil
    \geq \frac{n}{M} \geq \frac{\R}{2M}
    \left ( \frac{\R}{M} + 1 \right )
    \geq \frac{1}{2}\left \lfloor \frac{\R}{M} \right \rfloor
    \left ( \left \lfloor \frac{\R}{M} \right \rfloor
    + 1 \right ) = n_{\lfloor \R/M \rfloor}.
$$
Therefore, Theorem \ref{thm:superDiag} applies when $M = \lceil
\R^2 / (2n-\R) \rceil$, and we have
\begin{align*}
\qquad K^+(n,n-\R) & \ \leq \ K^+\left (\lceil n/M \rceil, \lceil
n/M \rceil - \lfloor
\R/M \rfloor \right )^M \\
& \ = \ \left ( \left \lfloor \frac{\R}{M} \right \rfloor + 1
\right )^M \leq \ \left ( 2n/\R \right )^{\lceil \R^2/(2n-\R)
\rceil}. \qquad \qquad \qquad \square
\end{align*}

We get the following corollary by letting $\R = (1-\lambda)n$.

\begin{cor} For any $\lambda$ with $0 \leq \lambda < 1$ and
$\lambda n$ integral,
$$
K^+(n,\lambda n) \leq \left ( \frac{2}{1-\lambda} \right )^{\lceil
n(1-\lambda)^2/(1+\lambda) \rceil}.
$$
\end{cor}

For each $\lambda$, this gives an exponential upper bound on $K^+(n,
\lambda n)$.  For example, when $\lambda = 1/2$, we have $K^+(n,n/2) \leq
4^{\lceil n/6 \rceil} < 1.26^{n+5}$.

\section{Difference bounds} \label{sec:difference}

In the discussion of asymptotic behavior above, we repeatedly used
the fact that $K^+(n,R)$ increases as $n$ increases or $R$
decreases. In fact, a cursory examination of the table included at
the end of this paper reveals that, at least above the diagonal,
the increase from entry to adjacent entry (i.e., to the right or
upward) is strict, and grows with increasing $n$ and $R$. Here, we
examine these ``difference'' patterns in more detail, by
considering the number of 0's or 1's in minimal codes.

\begin{prop} \label{prop:countingZeroes}
Let $\phi(n,R)$ be the maximum total number of 0's in a minimal
$(n,R)^+$-code, and let $\bar{\phi}(n,R)$ be the minimum number of
1's in a minimal $(n,R)^+$-code.  Furthermore, assume that $R \leq
n$.  Then we have the following: \smallskip

\indent 1. $K^+(n,R)-K^+(n-1,R) \geq \phi(n,R)/n $,

\indent 2. $ K^+(n-1,R) \leq \bar{\phi}(n,R)/n$,

\indent 3. $ K^+(n,R) < K^+(n+1,R)$, and

\indent 4. $ K^+(n,R) > K^+(n,R+1)$.
\end{prop}
\proof  \textbf{[1]} The proof is similar to that of Theorem
\ref{thm:superDiag}.  Let $\C$ be a minimal $(n,R)^+$-code.  Then
the codewords of $\C$ contain at least $\phi(n,R)$ 0's, and we may
choose a coordinate $i$ at which at least $\phi(n,R)/n$ codewords
of $\C$ have a 0.  The contraction $\C^i$ has at most
$K^+(n,R)-\phi(n,R)/n$ codewords, and downward $R$-covers
$Q_{n-1}$ by Lemma \ref{lem:RestrictCodes}.  [1] follows since
$K^+(n-1,R)\leq |\C^i|$.

\textbf{[2]} Let $\mathcal{C}$ be a minimal $(n,R)$-code achieving
$\bar{\phi}(n,R)$, and thus achieving $\phi(n,R)$, since
$\bar{\phi}(n,R)=nK^+(n,R)-\phi(n,R)$. Let $a_l$ be the number of
codewords of $\mathcal{C}$ at level $l$.  Then by part 1 we have
\begin{eqnarray}
K^+(n,R)-K^+(n-1,R) & \geq &  \frac{1}{n}
    \sum_{l=0}^{n}(n-l)a_l  \nonumber \\
& = & K^+(n,R) -  \frac{1}{n}\sum_{l=0}^n
    l\cdot a_l , \nonumber
\end{eqnarray}
and so
\begin{eqnarray}
K^+(n-1,R) & \leq &  \frac{1}{n}\sum_{l=0}^n
    l\cdot a_l  = \bar{\phi}(n,R)/n. \nonumber
\end{eqnarray}

\textbf{[3]} If $R \leq n$, then any $(n+1,R)^+$-code has at least two
codewords in it. It must therefore contain a vector other than
$\hat{1}$, so that
$$
\phi(n+1,R) \geq 1,
$$
and applying part 1 gives the desired result.

\textbf{[4]} Applying the direct sum bound
(\ref{eqn:directSumBound}),
\begin{eqnarray}
    \hskip 1in K^+(n,R+1) & \leq & K^+(n-1,R) \cdot K^+(1,1)  \\
    & = & K^+(n-1,R)    \nonumber \\
    & < & K^+(n,R).   \nonumber \hskip 2in \square
\end{eqnarray}

In order to get more out of Proposition
\ref{prop:countingZeroes}.1 than a difference of 1, we must more
carefully analyze the number of 0's in a code. A trivial lower
bound is obtained by noting that there are at most $\binom{n}{j}$
codewords with $j$ 0's.  A much better lower bound is obtained by
modifying the objective function in Proposition
\ref{prop:modifiedOneSideSphereCoverBound} to count the total
number of 0's in the code, as follows.

\begin{prop}\label{prop:phiLowerBound}  Let $IP_\phi^+(n,R)$ be
the result of the integer program in Proposition
\ref{prop:modifiedOneSideSphereCoverBound} with objective function
$\sum_{l=0}^n a_l$ replaced by $\sum_{l=0}^n (n-l)a_l$. Then for
any $n,R \geq 1$, $\phi(n,R) \geq IP_\phi^+(n,R)$.
\end{prop}

\section{Linear asymmetric codes} \label{sec:linear}

Up to this point, we have been considering general asymmetric
codes. However, a large portion of what is known about the
symmetric case concerns linear codes, so it is natural to ask what
can be said about asymmetric linear codes.  For example, for a
fixed radius, symmetric linear codes are asymptotically just as
efficient at covering the cube as nonlinear ones (up to a
multiplicative constant). The same statement is decidedly false,
however, in the asymmetric case. We will need some definitions
before we proceed with our results.

Let $\bar{\C}$, the \textit{1's complement} of $\C$, be the set
$\{\hat{1}-x\,|\,x \in \C\}$.  We say that an $(n,R)^+$-code $\C$
is a \textit{downward-asymmetric linear covering code of radius
$R$} (for short, an $[n,R]^+$-code) if it is a vector subspace of
$\mathbb{F}_2^n$, and $\C$ is an \textit{upward-asymmetric linear
covering code of radius $R$} (for short, an $[n,R]^-$-code) if its
1's complement is an $[n,R]^+$-code. Define $k^+[n,R]$ to be the
dimension of the smallest $[n,R]^+$-code, and $k^-[n,R]$ to be the
dimension of the smallest $[n,R]^-$-code.  In contrast with the
nonlinear case, we actually need to distinguish upward and
downward codes, as will become apparent shortly.

Call a code $\C$ \textit{self-complementary} if $\C = \bar{\C}$,
and define $k^\pm[n,R]$ to be the minimal dimension of a
self-complementary asymmetric linear code (for short, an
$[n,R]^\pm$-code). We do not need to specify ``upward'' or
``downward'' here, since a self-complementary code covers the cube
in one direction iff it covers it in the other. Finally, for a
code $\C$ and a coordinate $i$, define the {\em shortening} $\C_i
\subset Q_{n-1}$ of $\C$ to be the set of points of $\C$ with a
``0'' at position $i$, projected into $Q_{n-1}$ by deletion of
that bit. Thus $C_i$ is the 1's complement of the contraction at
$i$ of the 1's complement of $\C$.

We begin with a lemma which says that the dimension of a
downward-asymmetric linear covering code increases by at least $1$
when $n$ is incremented.

\begin{lemma}
$k^-[n-1,R] \leq k^-[n,R]-1$, for $n>R$.
\end{lemma}
\proof Let $\C$ be a minimal $[n,R]^-$-code.  Since $n>R$, one
$R$-ball cannot itself upward $R$-cover the entire $n$-cube.
Therefore $\C$ has at least two elements, and there exists a
coordinate $i$ where some vector $x \in \C$ has a 1.  The
shortening $\C_i$ of a linear code $\C$ is linear, and
$|\C_i|<|C|$ gives that $\dim(\C_i)<\dim(\C)$. \qed

The following theorem says exactly how large $k^+$, $k^-$, and
$k^\pm$ are.

\begin{theorem} For $n>0$,
$k^+[n,R] = k^\pm[n,R] = k^-[n,R] = \max\{1,n-R\}$.
\end{theorem}
\proof  We begin with the first equality.  Every $[n,R]^+$-code is
an $(n,R)^+$-code, so it must contain $\hat{1}$.  Containing
$\hat{1}$ is equivalent to being self-complementary for linear
codes, however, so every $[n,R]^+$-code is self-complementary.
That every $[n,R]^\pm$-code is an $[n,R]^+$-code is trivial, and
we have $k^+[n,R] = k^\pm[n,R]$.

Since every self-complementary asymmetric linear code is also an
upward-asymmetric linear code, we have $k^-[n,R] \leq k^\pm[n,R]$.
By induction using the previous lemma, we have $k^-[n,R] \geq
k^-[R,R]+n-R=n-R$ for $n > R$. Furthermore, when $n \leq R$,
$k^-[n,R]=0$.  Thus, for any $n>0$, $k^\pm[n,R] \geq k^-[n,R] \geq
\max\{1,n-R\}$, since all downward-directed codes include
$\hat{1}$.

To complete the proof, it suffices to find an $[n,R]^\pm$-code
$\mathcal{A}(n,R)$ of dimension $\max\{1,n-R\}$.  We construct one
inductively.  For $n\leq R+1$, let
$\mathcal{A}(n,R)=\{\hat{0},\hat{1}\}$.  For larger $n$, let
$\mathcal{A}(n,R)$ be the $[n,R]^{\pm}$-code
$\mathcal{A}(n-1,R)\oplus \{0,1\}$. \qed

\section{Conclusion} \label{sec:conclusion}

Open questions abound concerning asymmetric covering codes, since
the entirety of the theory of symmetric covering codes could be
reexamined in the asymmetric case.  However, several questions
stand out as particularly interesting.

It remains to determine the asymptotic order of magnitude of
$K^+(n,R)$ when neither $R$ nor $\R$ is constant (e.g., $R$ is
linear in $n$). Additionally, the analysis above concerned only
the case of binary codes. It seems a natural next step to
investigate asymmetric covering codes on more symbols than two.
Perhaps an interesting way to define the notion of $R$-balls in
that case would be to take the set of vectors in $\mathbb{Z}_a^n$
which differ from a given vector by \textit{increasing} at most
$R$ of its coordinates.  Or, maybe what should be asked for is
that the sum of the ``increases'' in each coordinate add up to at
most $R$. In fact, one could imagine asking similar questions of
much more general classes of posets: geometric lattices, Cartesian
products of some base set of posets, etc.  Furthermore, the
definition of $R$-ball used here certainly is not the only
imaginable one. Perhaps it would interesting to look at sets which
permit $R_1$ changes from 0's to 1's and $R_2$ changes from 1's to
0's, or which permit $\alpha\!\cdot\! w(x)$ (asymmetric or
symmetric) changes to $x$ for some $\alpha \in [0,1)$. Clearly
there is a lot of room for generalization.

The questions about classification of codes that arise in the
context of symmetric codes are relevant here as well.  What
possible forms do $[n,R]^\pm$-codes take?  Can anything be said
about how \textit{many} minimal $(n,R)^+$-codes there are? How
close to perfect -- i.e., no overlap between $R$-balls centered at
codewords -- can an $(n,R)^+$-code be?  Since linear codes, which
provide easily computable examples of efficient symmetric covering
codes, are so far from the sphere-covering lower bound in the
asymmetric case, does there exist a family of asymmetric codes
which are polynomial-time computable (in $n$) and which are within
a constant of optimal?

One important concrete problem stands: to find better upper and
lower bounds on $K^+(n,R)$ for small $n$ and $R$.  Table 1
demonstrates our best known bounds, to be interpreted as follows.
All entries weakly to the left or weakly below a subscript of `d'
are determined by Theorem \ref{thm:superDiag}. A subscript of `i'
denotes a lower bound found by the integer program in Proposition
\ref{prop:modifiedOneSideSphereCoverBound}. An `e' means an
explicit code was found exhaustively.  The subscript `m' indicates
that Proposition \ref{prop:phiLowerBound}, with the integer
program modified by an extra combinatorial constraint, was used to
compute a lower bound for $\phi(n,R)$, and then a lower bound for
$K^+(n,R)$ was found with Proposition \ref{prop:countingZeroes}.1.
An `s' means the upper bound is from a direct sum of codes of the
type in (\ref{eqn:directSumBound}). In particular, $K^+(10,4)\leq
K^+(5,2)\cdot K^+(5,2)$. Otherwise, no subscript on the left means
the lower bound was found using Proposition
\ref{prop:countingZeroes}.1 with $\phi(n,R)$ bounded below by
Proposition \ref{prop:phiLowerBound}, and no subscript on the
right means the upper bound corresponds to a code found greedily.
We note that Applegate, Rains and Sloane (\cite{ARS02}) already
claim the improvements $K^+(7,1)=31$, $K^+(8,1)=58$, $K^+(9,1)\leq
106$, $K^+(10,1)\leq 196$, $K^+(11,1)\leq 352$ and $K^+(12,1)\leq
670$.

\def\baselinestretch{1}
\begin{table}[!p]
\begin{center}
\begin{tabular}{|c||c|c|c|c|c|c|c|} \hline
$R \backslash n$ & 2 & 3 & 4 & 5 & 6 & 7  \\
\hline \hline 1  & 2 & 3$_\mathrm{d}$ &
$_\mathrm{i}$6$_\mathrm{e}$ & 10$_\mathrm{e}$ &
$_\mathrm{m}$18$_\mathrm{e}$& 30-34\\
\hline 2   & 1 &  2 &  3$_\mathrm{d}$ &  5$_\mathrm{e}$ &
$_\mathrm{m}$8$_\mathrm{e}$ & 13-15$_\mathrm{e}$  \\
\hline 3   & 1  & 1 &  2 &  3 &  4$_\mathrm{d}$  &
$_\mathrm{i}$6-7  \\
\hline 4   & 1  & 1  & 1  & 2  & 3 &  4$_\mathrm{d}$  \\
\hline 5   & 1  & 1  & 1 &  1 &  2  & 3  \\
\hline 6& 1 &1 &  1  & 1  & 1 &  2 \\
\hline \multicolumn{7}{c}{} \\\hline $R \backslash n$ &  8 & 9 &
10 & 11 &
12 & 13 \\
\hline \hline 1  &  52-67 & 93-121
& 162-229 & 306-433 & 563-813 & 1046-1626$_\mathrm{s}$ \\
\hline 2   &   20-25 & 32-46 &
52-81 & 87-141 & 148-262 & 254-524$_\mathrm{s}$ \\
\hline 3   &  $_\mathrm{i}$9-13  & $_\mathrm{i}$14-21 &
22-36 & 34-64 &54-105 & 88-210 \\
\hline 4   & 6 &8-11 & 12-16$_\mathrm{s}$
&17-30 &26-49 &40-83 \\
\hline 5   &  4$_\mathrm{d}$ &6 &8-9& 11-16&
15-27 &22-48 \\
\hline 6&  3  & 4 &5$_\mathrm{d}$ &7-8& 10-15&
14-23  \\
\hline 7 &2 &3 &  4 &  5$_\mathrm{d}$ &7$_\mathrm{e}$ &9-12\\
\hline 8  & 1& 2 &3& 4& 5$_\mathrm{d}$ &7 \\
\hline 9 &1& 1& 2& 3& 4& 5$_\mathrm{d}$ \\
\hline 10 & 1& 1 &1 &2& 3 &4 \\
\hline 11 &1 &1& 1 &1 &2& 3 \\ \hline
\end{tabular}\caption{Best known bounds for $K^+(n,R)$}
\end{center}\label{tab:LPFirst}
\end{table}

\end{document}